\magnification=\magstep1
\input amstex
\documentstyle{amsppt}

\define\defeq{\overset{\text{def}}\to=}

\def \isom {\buildrel \sim \over \rightarrow}

\def \Im{\operatorname {Im}}
\def \Pic{\operatorname {Pic}}
\def \pr{\operatorname {pr}}
\def \sat{\operatorname {sat}}
\def \tor{\operatorname {tor}}
\def \new{\operatorname {new}}

\def \Sp{\operatorname {Sp}}

\def \ord{\operatorname {ord}}
\def \Ker{\operatorname {Ker}}
\def \char{\operatorname {char}}
\def \Spec{\operatorname {Spec}}
\def \acts\ trivially\ on{\operatorname {acts\ trivially\ on}}
\def \closed\ subgroups{\operatorname {\closed\ subgroups}}
\def \Sp{\operatorname {Sp}}
\def \sp{\operatorname {Sp}}

\NoRunningHeads
\NoBlackBoxes
\topmatter

\title
Variation of fundamental groups of curves in positive characteristic%s
\endtitle
\bigskip

\author
{MOHAMED SA\"IDI and AKIO TAMAGAWA}
\endauthor

\abstract In this paper we prove that given a non-isotrivial family of 
%proper and smooth 
hyperbolic curves
%in characteristic $p>0$, then 
in positive characteristic, 
the isomorphism type of the geometric fundamental group is not constant on the fibres of the family.
\endabstract

\toc

\subhead
\S0. Introduction
\endsubhead

\subhead
\S1. Review of the Sheaf of Locally Exact Differentials in Characteristic $p>0$ and its Theta Divisor 
\endsubhead

\subhead
\S2. The Specialisation Theorem for Fundamental Groups
\endsubhead

\subhead
\S3. Families of Curves with a Given Fundamental Group in Positive Characteristic 
\endsubhead

\subhead
\S4. Proof of the Main Theorems
\endsubhead

\endtoc

\endtopmatter
\document

\subhead
\S 0. Introduction
\endsubhead
Let $k$ be an algebraically closed field and 
$X$ a proper, smooth, and connected 
%hyperbolic 
algebraic curve over $k$ of genus $g\ge 2$.
The structure of the \'etale fundamental group $\pi_1(X)$ of $X$ is well understood, 
if $\char(k)=0$, thanks to Riemann's Existence Theorem.
Namely, $\pi_1(X)$ is isomorphic to the profinite completion $\Gamma _g$ of the topological fundamental group of a compact 
orientable topological surface of genus $g$.
In particular, the 
%structure 
isomorphism type 
of $\pi_1(X)$ is constant and depends only on $g$ in this case.  In the case where $\char(k)=p>0$, the structure of the full $\pi_1(X)$ is still 
mysterious and is far from being understood. One only knows the structure of certain quotients of $\pi_1(X)$ in this case. More precisely, let $\pi_1(X)^p$ (resp. $\pi_1(X)^{p'}$)
be the maximal pro-$p$ (resp. maximal pro-prime-to-$p$) quotient of $\pi_1(X)$. Then it is well-known that there exists a surjective continuous specialisation homomorphism
$\Sp: \Gamma _g\twoheadrightarrow \pi_1(X)$, which induces an isomorphism $\Sp': \Gamma _g^{p'}\isom {\pi_1(X)}^{p'}$
between the maximal pro-prime-to-$p$ parts (cf. [SGA1]). Moreover, $\pi_1(X)^p$ is a free pro-$p$ group on 
$r
%=r_X
$-generators where 
%$r_X$ 
$r$ 
is the $p$-rank of 
(the Jacobian of) 
the curve $X$ (cf. [Sh]).
The full structure of $\pi_1(X)$ is not known for a single example of a curve $X$ of genus $g\ge 2$ in characteristic $p>0$.

In order to understand the complexity of the geometric fundamental group $\pi_1$ of hyperbolic curves in positive 
characteristic, it is natural to investigate the variation of the structure of $\pi_1$ 
when curves vary in their moduli. Let $\Cal M_{g,\Bbb F_p}$ be the coarse moduli space of 
proper, smooth, and connected curves of genus $g$ in characteristic $p>0$.
Given a point $x\in \Cal M_{g,\Bbb F_p}$, choose 
a geometric point $\bar x$ above $x$ and let $C_{\bar x}$ be a curve 
%over (the source of) $\bar x$ 
corresponding to the moduli point $\bar x$ (well-defined up to isomorphism). 
Then the isomorphism type of 
the (geometric) \'etale fundamental group $\pi_1(C_{\bar x})$ 
is independent of the choice of $\bar x$ and $C_{\bar x}$ (and the implicit 
base point on $C_{\bar x}$ used to define $\pi_1(C_{\bar x})$). 
(See [S], \S 4 for more details). 
An important tool in studying fundamental groups in positive characteristic is the specialisation theory of 
Grothendieck (cf. [SGA1]).
Given points $x,y\in \Cal M_{g,\Bbb F_p}$, such that $x\in \overline {\{y\}}$ holds, there exists a continuous surjective {\bf specialisation homomorphism} $\Sp_{y,x}:\pi_1(C_{\bar y})\to \pi_1(C_{\bar x})$. Concerning this specialisation homomorphism we have the following fundamental specialisation theorem (cf. Theorem 2.1), which is proven in [T3].

\proclaim {Theorem A}Let  $x,y\in \Cal M_{g,\Bbb F_p}$ be distinct points of $\Cal M_{g,\Bbb F_p}$ such that $x\in \overline {\{y\}}$ holds. Assume that $x$ is a {\bf closed} 
point of $\Cal M_{g,\Bbb F_p}$. Then the specialisation homomorphism $\Sp_{y,x}:\pi_1(C_{\bar y})\to \pi_1(C_{\bar x})$ is not an isomorphism.
\endproclaim

It is quite plausible that Theorem A may hold in general without the extra assumption that $x$ is a closed point of $\Cal M_{g,\Bbb F_p}$.
It is also plausible that the isomorphism type of $\pi_1(C_{\bar x})$ tends to depend on the moduli point $x\in \Cal M_{g,\Bbb F_p}$. 
In the spirit of Grothendieck's anabelian geometry this would suggest the possibility that strong anabelian phenomena for curves 
over algebraically closed fields of positive characteristic may occur, in contrast to the situation in characteristic $0$ where the geometric fundamental 
group carries only topological informations.

A slightly weaker approach to the above specialisation theorem is the following. Let $k$ be a field of characteristic $p>0$ and set $\Cal M_{g,k}\defeq \Cal M_{g,\Bbb F_p}
\times _{\Bbb F_p}k$. Let $\Cal S\subset \Cal M_{g,k}$ be a subvariety. 
We say that the (geometric) fundamental group $\pi_1$ is {\bf constant} on $\Cal S$ 
if, for any two points $y$ and $x$ of $\Cal S$, 
such that $x\in\overline{\{y\}}$ holds, the specialisation homomorphism 
$\Sp_{y,x}:\pi_1(C_{\bar y})\to \pi_1(C_{\bar x})$ is an isomorphism. 
%one has $\pi_1(C_{\bar x})\simeq\pi_1(C_{\bar y})$. 
We say that $\pi_1$ is not constant on $\Cal S$ if the contrary holds (cf. Definition 3.1). 
In [S] was raised the following question.

\definition {Question} Does $\Cal M_{g,k}$ contain any subvariety of positive dimension 
on which $\pi_1$ is constant?
\enddefinition

When $k$ is an algebraic closure of $\Bbb F_p$, the answer to this question is negative by Theorem A. 
For general $k$, the following result is proven in [S] (cf. [S], Theorem 4.4).

\proclaim{Theorem B} Let $k$ be a field of characteristic $p>0$ and $\Cal S\subset \Cal M_{g,k}$ a {\bf complete} 
subvariety of $\Cal M_{g,k}$ of positive dimension.
Then the (geometric) fundamental group $\pi_1$ is not constant on $\Cal S$.
\endproclaim

The aim of this paper is to remove the completeness assumption in Theorem B and give a negative answer to the above question in general.
The main result of this paper is the following (cf. Theorem 3.6).

\proclaim{Theorem C} Let $k$ be a field of characteristic $p>0$ and $\Cal S\subset \Cal M_{g,k}$ a 
(not necessarily complete) subvariety of $\Cal M_{g,k}$ of positive dimension. 
Then the (geometric) fundamental group $\pi_1$ is not constant on $\Cal S$.
\endproclaim

Note that the validity of Theorem A in general, i.e., without the extra assumption that $x$ is a closed point, would immediately imply Theorem C.
Our proof of Theorem C is quite different from the proof of Theorem B in [S], and relies on Theorem A 
and Raynaud's theory of theta divisors. We also prove certain variants of Theorem C for curves which 
are not necessarily proper  
(cf. Theorems 3.12 and 3.13). 

Next, we briefly review the contents of each section. In $\S1$ we review basic facts on the theta divisor of the sheaf of locally exact differentials on a curve in characteristic $p>0$.
In $\S2$ we review some key facts (proven in the course of proving Theorem A in [T3]) 
which are used in this paper. In $\S3$ we state the main theorems, 
and in $\S4$ we proceed to their proof.

\subhead
\S 1. Review of the Sheaf of Locally Exact Differentials in Characteristic $p>0$ and its Theta Divisor 
\endsubhead
In this paper $p$ denotes a (fixed) prime number. 
In this section, we will review Raynaud's theory of theta divisors in characteristic $p$,  
initiated in [R]. 

Let $S$ be an $\Bbb F_p$-scheme. 
We denote by $F_S$ the absolute Frobenius endomorphism $S@>>> S$. For an $S$-scheme $X$, we define $X_1$
to be the pull-back of $X$ by $F_S$. Thus, we have a cartesian square
$$
\CD
X_1@>>>  X\\
@VVV     @VVV \\
S @>{F_S}>> S\\
\endCD
$$
The absolute Frobenius endomorphism $F:X\to X$ induces in a natural way an $S$-morphism $F_{X/S}:X\to X_1$, the relative Frobenius 
morphism, which is an integral radicial morphism. 
Next, assume that $X$ is a proper and smooth 
%and connected 
(relative) 
{\bf $S$-curve} of genus $g$, i.e., the morphism $X\to S$ is 
proper and smooth and its 
fibres are (geometrically connected) curves of (constant) genus $g$. Then $X_1$ is also a proper and smooth 
$S$-curve of genus $g$, and $F_{X/S}:X\to X_1$ is finite locally free of degree $p$. 
The canonical differential 
$$(F_{X/S})_{*}d:(F_{X/S})_{*} \Cal O_X\to     (F_{X/S})_{*}\Omega ^1_{X/S}$$
is a morphism of $\Cal O_{X_1}$-modules. Its image 
$$B_X\defeq \Im ((F_{X/S})_{*}d)$$ 
is the sheaf of locally exact differentials. We have a natural exact sequence
$$0\to \Cal O_{X_1}\to (F_{X/S})_{*} \Cal O_X\to B_X\to 0,$$
and $B_X$ is a locally free $\Cal O_{X_1}$-module of rank $p-1$.

Let $J$ (resp. $J_1$) denote the relative Jacobian $\Pic^0_{X/S}$ of $X$ (resp. $\Pic^0_{X_1/S}$ of $X_1$) over $S$,
which is an abelian scheme of relative dimension $g$ over $S$ (cf. [BLR], 9.4, Proposition 4). Then, \'etale-locally on $S$, there exists a universal 
degree $0$ line bundle  $\Cal L_1$ on $X_1\times _SJ_1$. We define $\Theta= \Theta _X$ to be the closed subscheme of $J_1$ defined
by the $0$-th Fitting ideal of $R^1(\pr_{J_1})_{*}(\pr_{X_1}^*(B_X)\otimes\Cal L_1)$, where $\pr_{X_1}$ and $\pr_{J_1}$ denote the projections 
$X_1\times _SJ_1\to X_1$ and $X_1\times _SJ_1\to J_1$, respectively. 
Since $\Theta$ is independent of the choice of $\Cal L_1$ 
(cf. [T2], Proposition 2.2(i)), we can define $\Theta$ not only \'etale-locally but also globally on $S$. 
By definition, the formation of $\Theta$
commutes with any base change of $S$. The following result is essentially due to Raynaud.

\proclaim {Theorem 1.1} $\Theta$ is a relative Cartier divisor on $J_1/S$.
\endproclaim

\demo {Proof} See [T3], Theorem (5.1).
\qed
\enddemo

In this paper we will refer to the divisor $\Theta$ as the (relative) {\bf Raynaud theta divisor}.

{}From now on, we assume $S=\Spec (k)$, where $k$ is an algebraically closed field of characteristic $p>0$. 

\definition {Definition 1.2} Let $M$ be an abelian group and $M_{\tor}$ the subgroup of torsion elements of $M$.

(i)\ For each element $x\in M_{\tor}$ we define the subset $\sat (x)$ (which we call the 
{\bf saturation} of $x$) of $M$ to be the set
of elements in the form $i\cdot x$, where $i$ is an integer prime to the order $N_x$ of $x$. 
(Thus, $\sharp (\sat(x))=\varphi(N_x)$.) 

(ii) For each subset $X$ of $M_{\tor}$ we define the subset $\sat(X)$ (which we call the {\bf saturation} 
of $X$) of $M$ to be the union of $\sat(x)$ for all $x\in X$. 
(Note that $\sat(X)\subset M_{\tor}$.) 
Moreover, we say that $X$ is saturated if $\sat(X)=X$.
\enddefinition

Recall that given a scheme $S$, a cyclic group $G$ of order $N$ which is invertible on $S$, and an abelian $S$-scheme $\Cal A$ endowed with a $G$-action,
one defines naturally the {\bf ``new part''} $\Cal A_{\new}$ of $\Cal A$ with respect to this action (cf. [T3], 
$\S4$ for more details).  
The following result relates the geometry of the 
Raynaud theta divisor to fundamental groups.

\proclaim {Proposition 1.3} Let $N$ be a positive integer prime to $p$. Let $x$ be a torsion point of $J_1(k)$ of order $N$, and $Y_1\to X_1$ the $\mu_N$-torsor
associated to $x$. Then, $\sat (x)\cap \Theta(k)=\varnothing$ if and only if $Y_1\to X_1$ is new-ordinary in the sense that $(J_{Y_1})_{\new}$ is an ordinary abelian variety.
\endproclaim 

\demo{Proof} See [T3], Proposition (5.2).
\qed
\enddemo
 
 % We also have the following result concerning the geometry of the Raynaud theta divisor $\Theta$.
% \proclaim {Proposition 1.4} Assume $g>0$.
% 
% (i)\ $\Theta$ is algebraically equivalent to $(p-1)\Theta _{1,\cl}$, where $\Theta _{1,\cl}$ is 
%a classical theta divisor for $X_1$.
% 
% (ii)\ $\Theta$ is an ample divisor.
% 
% (iii)\ $(X_1,\Theta)=(p-1)g$, where $X_1$ is regarded as a subvariety of $J_1$ by means of an albanese embedding.
% \endproclaim
%
%\demo{Proof} See [T3], Proposition (5.5).
%\qed
%\enddemo

\subhead
2. The Specialisation Theorem for Fundamental Groups
\endsubhead
In this section $k_0$ denotes an algebraic closure of the prime field $\Bbb F_p$ of characteristic $p>0$.
Let $S$ be an $\Bbb F_p$-scheme, $s$ and $t$ points of $S$ such that $s\in \overline{\{t\}}$ holds.
We denote by $\bar s$ and $\bar t$ geometric points above $s$ and $t$, respectively. Let $X$ be a proper and smooth 
%(relative) 
$S$-curve of genus $g$.
Write $X_{\bar s}\defeq X\times _S\bar s$, $X_{\bar t}\defeq X\times _S\bar t$ 
for the geometric fibres of $X$ above $s$ and $t$, respectively.
Then we have a {\bf specialisation homomorphism} (cf. [SGA1], Expos\'e IX, 4 and Expos\'e XIII, 2.10)
$$\Sp:\pi_1(X_{\bar t})\to \pi_1(X_{\bar s}),$$
which is surjective and induces an isomorphism 
$\pi_1(X_{\bar t})^{(p')}\isom \pi_1(X_{\bar s})^{(p')}$ 
between the maximal pro-prime-to-$p$ quotients.

Recall that a curve over a field containing $k_0=\overline {\Bbb F}_p$ is {\bf constant}, if it descends to a curve over $k_0$. The following result %proven in [T3]
is fundamental.

\proclaim{Theorem 2.1} Assume $g\ge 2$. Assume that $X_{\bar s}$ is constant and that $X_{\bar t}$ is not constant. 
Then the specialisation homomorphism 
$\Sp:\pi_1(X_{\bar t})\to \pi_1(X_{\bar s})$ is not an isomorphism.
\endproclaim

\demo{Proof} See [T3], Theorem (8.1). 
%, where a tame version of the above statement is also proven.
\qed
\enddemo

Theorem 2.1 follows easily from the following (cf. [T3], Theorem (6.1)).

\proclaim{Theorem 2.2} Let $R$ be a complete discrete valuation ring isomorphic to $k_0[[t]]$, and 
set $S\defeq\Spec(R)=\{\eta,s\}$, where
$\eta$ (resp. $s$) stands for the generic (resp. closed) point of $S$. Let $X$ be a proper and smooth 
%(relative) 
$S$-curve of genus $g\ge 2$,
and assume that $X_{\bar \eta}$ is not constant. Then the specialisation homomorphism 
$\Sp:\pi_1(X_{\bar \eta})\to \pi_1(X_{\bar s})$ is not an isomorphism.
\endproclaim

For an abelian variety $A$ over an algebraically closed field $k$, 
%of characteristic $p$, 
we write $A\{p'\}$ for the union of $A[N](k)$ for all positive integers $N$ prime to $p$,
and for a subscheme $Z$ of $A$ we write $Z\{p'\}\defeq Z(k)\cap A\{p'\}$.

Let $X\to S=\Spec(R)=\{\eta,s\}$ be as in Theorem 2.2. 
(In the following discussion including Lemma 2.3, 
$R$ may be an arbitrary discrete valuation 
ring of equal characteristic $p>0$.)  
%Suppose that the specialisation homomorphism  
%$$\Sp:\pi_1(X_{\bar \eta})\to \pi_1(X_{\bar s})$$ 
%is an isomorphism. Then, since a radicial morphism does not change the fundamental group, we see that 
%$$\Sp_1:\pi_1(X_{1,\bar \eta})\to \pi_1(X_{1,\bar s})$$ 
%is also an isomorphism, where $X_1$ denotes the Frobenius twist of $X$ over $S$ (cf. $\S1$). 
Let $X_1$ be the Frobenius twist of $X$ over $S$ (cf. $\S1$) and 
$J_1$ the (relative) Jacobian of $X_1$ over $S$. This is an abelian scheme over $S$ and can be identified with the N\'eron model of $J_{1,\eta}\defeq J_1\times _S{\eta}$.
By the N\'eron property, we have a natural specialisation isomorphism 
$$J_{1,\bar \eta}\{p'\}\to J_{1,\bar s}\{p'\},$$
where $J_{1,\bar \eta}\defeq J_1\times _S{\bar \eta}$ and 
$J_{1,\bar s}\defeq J_1\times _S{\bar s}$. 
Identifying these two abelian groups with each other by this specialisation isomorphism, we will write
$$J_1\{p'\}\defeq J_{1,\bar \eta}\{p'\}=J_{1,\bar s}\{p'\}.$$
(Thus, $J_1\{p'\}$ is a mere abelian group). Moreover, we have the Raynaud theta divisor $\Theta$ in $J_1$, and under the above identification we have
$$\Theta _{\bar \eta}\{p'\}\subset \Theta _{\bar s}\{p'\}\ ({\text {in}}\ J_1\{p'\}).$$
The following is a crucial observation (cf. [T3], Lemma (6.2)). 

\proclaim {Lemma 2.3}
If $\Sp$ is an isomorphism, then $\sat (\Theta _{\bar \eta}\{p'\})=\sat (\Theta _{\bar s}\{p'\})$ 
must hold in $J_1\{p'\}$.
\endproclaim

In the course of proving Theorem 2.2, the following more precise statement is proven (cf. [T3], $\S 7$).

\proclaim{Proposition 2.4} 
Under the assumptions of Theorem 2.2 ($R\simeq k_0[[t]]$), 
there exists 
%an open subgroup $H_{\bar \eta}'\subset \pi_1(X_{1,\bar \eta})^{(p')}$, 
%corresponding to an open subgroup 
%$H_{\bar s}'\subset \pi_1(X_{1,\bar s})^{(p')}$, $H_{\bar \eta}\subset \pi_1(X_{1,\bar \eta})$ 
%the inverse image of $H_{\bar \eta}'$ in $\pi_1(X_{1,\bar \eta})$ , 
%corresponding to the open subgroup $H_{\bar s}\subset \pi_1(X_{1,\bar s})$, $H_{\bar \eta}$ 
%corresponds to 
a finite \'etale cover $Y\to X$ whose Galois closure is of degree prime to $p$,  
such that 
$$\sat (\Theta _{Y,\bar \eta}\{p'\})\subsetneq \sat (\Theta _{Y,\bar s}\{p'\})$$
holds in $J_{Y_1}\{p'\}$. 
Here, $J_{Y_1}$ and $\Theta _{Y}
\subset J_{Y_1}
$ denote the (relative) Jacobian of $Y_1$ over $S$ 
and the (relative) Raynaud theta divisor 
for $Y\to S$, respectively. 
\endproclaim

%\proclaim{Proposition 2.*} 
%Let $R_0$ be a complete discrete valuation ring isomorphic to 
%$\overline{\Bbb F}_p[[t]]$, set $S_0\defeq\Spec(R_0)$, and 
%let $f_0: X_0\to S_0$ be a proper smooth 
%%(relative) 
%$S_0$-curve 
%of (constant) genus $g\geq 2$. Let $\eta_0$ and $s_0$ be the 
%generic and the closed points of $S_0$, respectively. 
%Assume that $(X_0)_{\bar \eta_0}$ does not descend to 
%%a curve over 
%$\overline{\Bbb F}_p$. 
%Then there exists 
%a connected finite \'etale cover $Y_0\to X_0$ 
%whose Galois closure is geometrically connected over $S_0$ and 
%of degree prime to $p$ over $X_0$, 
%such that the $p$-rank of 
%%the geometric fibre 
%$(Y_0)_{\bar s_0}$ 
%%of $Y_0$ above $s_0$ 
%is (strictly) smaller than 
%the $p$-rank of 
%%the geometric fibre 
%$(Y_0)_{\bar \eta_0}$. 
%%of $Y_0$ above $\eta_0$. 
%\endproclaim

\demo{Remark 2.5}
In fact, the finite \'etale cover $Y\to X$ in Proposition 2.4 can be chosen to factorise 
as $Y=X_3\to X_2\to X_1\to X_0=X$, where $X_i\to X_{i-1}$ is a $\mu_{N_i}$-torsor
for a suitable positive integer $N_i$ prime to $p$, $i\in \{1,2,3\}$. 
\enddemo

\subhead
\S 3. Families of Curves with a Given Fundamental Group in Positive Characteristic 
\endsubhead
In this section we state our main results. We will start with the following elementary definition. 

\definition{Definition 3.1} 
Let $S$ be a set and assume that for each $s\in S$, a profinite group $\Pi_s$ is given. 
We denote by $\Pi$ the map from $S$ to the set of isomorphism classes of profinite groups 
that assigns to each $s\in S$ the isomorphism class of $\Pi_s$. 

(i) We say that $\Pi$ is {\bf constant} on $S$ if, for any $s,t\in S$, one has $\Pi_s\simeq\Pi_t$. 

(ii) Assume moreover that $S$ is the underlying (vertex) set 
%(of vertices) 
of an oriented graph, and 
that for each oriented 
edge $t\rightarrowtail s$ linking vertices $t$ and $s$ 
of the graph $S$, a surjective (continuous) homomorphism 
$\sp_{t,s}: \Pi_t\to \Pi_s$ is given. Then we say that $\Pi$ is {\bf $\sp$-constant} on $S$ 
if for any oriented edge $t\rightarrowtail s$ of the graph $S$, the homomorphism 
$\sp_{t,s}: \Pi_t\to\Pi_s$ is an isomorphism. 
\enddefinition

\proclaim{Lemma 3.2} In the situation of Definition 3.1 (ii), consider the following 
conditions: 
(i) $\Pi$ is constant; and 
(ii) $\Pi$ is $\sp$-constant. 
If $S$ is connected (as a graph), then one has (ii)$\implies$(i). 
If $\Pi_s$ is finitely generated (as a profinite group) for each $s\in S$, 
then one has (i)$\implies$(ii). 
\endproclaim

\demo{Proof} 
The first assertion is clear. The second assertion follows from the fact that 
a finitely generated profinite group is Hopfian (cf. [FJ], Proposition 15.3). 
\qed\enddemo

In the rest of this section, let $p$ be a prime number and 
$g$ an integer $\geq 2$ (unless otherwise stated), and set 
$\Cal M_g\defeq \Cal M_{g,\Bbb F_p}
%\to \Bbb F_p
$, the {\bf coarse moduli space} of proper, smooth and geometrically connected curves of genus $g$ 
in characteristic $p
%>0
$.
Let $k$ be a field of characteristic 
$p
%>0
$, and $\Cal M_{g,k}\defeq \Cal M_g\times _{\Bbb F_p}k$, which turns out to be the 
coarse moduli space of proper, smooth and geometrically connected curves of genus $g$ 
over $k$-schemes.
%Recall that if $x \in \Cal M_{g,k}$ is a point 
%then one associates naturally to $x$ a geometric fundamental group 
%$\pi_1(C_{\bar x})$ 
Given a point $x\in \Cal M_{g,k}$, choose 
a geometric point $\bar x$ above $x$ and let $C_{\bar x}$ be a curve 
%over (the source of) $\bar x$ 
corresponding to the moduli point $\bar x$ (well-defined up to isomorphism). 
Then the isomorphism type of the (geometric) \'etale fundamental group $(\pi_1)_x\defeq \pi_1(C_{\bar x})$, 
which is a finitely generated profinite group 
(cf. [SGA1], Expos\'e X, Th\'eor\`eme 2.6),  
is independent of the choice of $\bar x$ and $C_{\bar x}$ (and the implicit 
base point on $C_{\bar x}$ used to define $\pi_1(C_{\bar x})$). 
If $x$ and $y$ are points in $\Cal M_{g,k}$ such that $x\in \overline {\{y\}}$ holds, 
then Grothendieck's specialisation theory for fundamental groups
implies the existence of a (continuous) surjective specialisation homomorphism 
$%\Sp=
\Sp_{y,x}:\pi_1(C_{\bar y})\twoheadrightarrow \pi_1(C_{\bar x})$. 
(See [S], $\S 4$ for more details.) 

Similarly, 
let $S$ be a scheme of characteristic $p$ and $f:X\to S$ 
a proper and smooth 
%(relative) 
$S$-curve of genus $g$. 
Given a point $s\in S$, choose 
a geometric point $\bar s$ above $s$ and let $X_{\bar s}$ be the 
geometric fibre of $f$ at $\bar s$. 
Then the isomorphism type of 
the (geometric) \'etale fundamental group $(\pi_1)_s\defeq \pi_1(X_{\bar s})$, 
which is a finitely generated profinite group,  
is independent of the choice of $\bar s$ (and the implicit 
base point on $X_{\bar s}$ used to define $\pi_1(X_{\bar s})$). 
If $s$ and $t$ are points in $S$ such that $s\in \overline {\{t\}}$ holds, 
then Grothendieck's specialisation theory for fundamental groups
implies the existence of a (continuous) surjective specialisation homomorphism 
$%\Sp=
\Sp_{t,s}:\pi_1(X_{\bar t})\twoheadrightarrow \pi_1(X_{\bar s})$. 

%\definition {Definition 3.1} 
%(i)\ Let $\Cal S\subset \Cal M_{g,k}$ be a non-empty subscheme of $\Cal M_{g,k}$. 
%We say that the (geometric) fundamental group
%$\pi_1$ is {\bf constant} on $\Cal S$ if for any two points $x$ and $y$ of $\Cal S$ 
%such that $x\in \overline {\{y\}}$ holds, then the corresponding specialisation homomorphism
%$\Sp_{y,x}:\pi_1(C_{\bar y})\to \pi_1(C_{\bar x})$ is an isomorphism. We say that $\pi_1$ is 
%{\bf not constant} on $\Cal S$ if the contrary holds namely: there exist
%two points $x$ and $y$ of $S$, such that $x\in \overline {\{y\}}$ holds and that 
%the corresponding specialisation homomorphism $\Sp_{y,x}:\pi_1(C_{\bar y})\to \pi_1(C_{\bar x})$ is 
%not an isomorphism. 
%
%(ii)\ Let $S$ be a connected scheme of characteristic $p$, and $f:X\to S$ a 
%%(relative) 
%proper and smooth 
%$S$-curve of genus $g$. We say that the (geometric) 
%fundamental group is {\bf constant} on the family $f$ if for any two points $s$ and $t$ of $S$, 
%such that $s\in \overline {\{t\}}$ holds, the corresponding 
%specialisation homomorphism $\Sp_{t,s}:\pi_1(X_{\bar t})\to \pi_1(X_{\bar s})$ is an isomorphism, 
%where $X_{\bar t}\defeq X\times _S \overline {k(t)}$
%and $X_{\bar s}\defeq X\times _S \overline {k(s)}$ are the geometric fibres of $f$ above the points $t$ 
%and $s$, respectively. Otherwise, i.e., if this condition is not 
%satisfied, we say that the fundamental group is {\bf not constant} on $f$.
%\enddefinition

Note that any topological space $T$ 
(e.g., any subset of a scheme) 
can be regarded as an oriented graph by 
assigning 
a vertex to each point of $T$ and 
an oriented edge $t\rightarrowtail s$ to each pair $(s,t)$ of distinct points of $T$ 
such that $s\in\overline{\{t\}}$ holds. 

%Note that if $S$ is a connected scheme of characteristic $p$ and 
%$f:X\to S$ a 
%%(relative) 
%proper and smooth $S$-curve of genus $g$, then the (geometric) fundamental group
%is constant on the family $f$ if and only if the (geometric) fundamental group 
%is constant on the image of the classifying map $S\to \Cal M_g$
%of the family $f$. 

\proclaim{Lemma 3.3}
(i) Let $\Cal S$ be a connected subscheme of $\Cal M_{g,k}$. Then $\pi_1$ is constant on $\Cal S$ 
if and only if $\pi_1$ is $\Sp$-constant on $\Cal S$ (cf. Definition 3.1). 

(ii) Let $k$ be a field of characteristic $p$, 
$S$ a connected $k$-scheme of finite type, 
$f:X\to S$ a proper and smooth 
%(relative) 
$S$-curve of genus $g$ 
and $h:S\to \Cal M_{g,k}$ the (coarse) classifying morphism for $f$. 
Then the following are all equivalent: 

(a) $\pi_1$ is constant on $S$; 

(b) $\pi_1$ is $\Sp$-constant on $S$; %and 

(c) $\pi_1$ is constant on $h(S)$. 
\endproclaim

\demo{Proof} This follows from Lemma 3.2, together with the (easily verified) fact 
that a noetherian scheme is connected as a scheme 
(or, equivalently, as a topological space) 
if and only if 
%it is connected as a graph. 
the associated graph is connected. 
\qed\enddemo

In [S] was raised the following question.

\definition {Question 3.4} Let $k$ be a field of characteristic $p$. 
Does $\Cal M_{g,k}$ contain any subvariety of positive dimension on which 
%the (geometric) fundamental group 
$\pi_1$ is constant?
\enddefinition

The following result is proven in [S] (cf. [S], Theorem 4.4.). 

\proclaim{Theorem 3.5} 
Let $k$ be a field of characteristic $p$ and $\Cal S\subset \Cal M_{g,k}$ a {\bf complete} subvariety 
of $\Cal M_{g,k}$ of positive dimension. 
Then 
%the (geometric) fundamental group 
$\pi_1$ 
is not constant on $\Cal S$.
\endproclaim

The main result of this paper is the following theorem, 
where we remove the assumption in Theorem 3.5 that the subvariety $\Cal S$ is complete.

\proclaim{Theorem 3.6} Let $k$ be a field of characteristic $p$, and $\Cal S\subset \Cal M_{g,k}$ a 
(not necessarily complete and even not necessarily closed) subvariety of $\Cal M_{g,k}$ of positive
dimension. Then 
%the (geometric) fundamental group 
$\pi_1$ is not constant on $\Cal S$.
\endproclaim

%As in the proof of Theorem 3.5 (cf. [S], Proof of Theorem 4.4) 
Theorem 3.6 follows 
%easily 
{}from Theorem 3.9 below. 

In the rest of this section, 
let $k$ be a field of characteristic $p$, 
$S$ a connected and reduced $k$-scheme of finite type,
and 
$f:X\to S$ a proper and smooth 
%(relative) 
$S$-curve of (constant) genus $g$. 

\definition{Definition 3.7} We say that $f$ is {\bf isotrivial}, if there exist 
a finite extension $k'/k$, a connected $k'$-scheme $S'$, a finite \'etale $k'$-morphism 
$S'\to S\times_kk'$, and a proper and smooth $k'$-curve $X_0'$, such that 
$X\times_S S'$ is isomorphic to $X_0'\times_{k'}S'$ over $S'$. 
\enddefinition

\proclaim{Lemma 3.8} 
%Assume that $g\geq 2$ and write 
Write $h: S\to\Cal M_{g,k}$ 
for the (coarse) classifying morphism for $f$. 
%and consider the following conditions: 
%
%(i) $f$ is isotrivial; and 
%
%(ii) the (coarse) classifying morphism $h: S\to\Cal M_{g,k}$ for $f$ is set-theoretically constant. 
%
%\noindent
%Then one has (i)$\implies$(ii). If, moreover, $S$ is reduced, then one has (ii)$\implies$(i). 
Then the following conditions are all equivalent. 

(i) $f$ is isotrivial. 

(ii) The image of $h$ consists of a single closed point of $\Cal M_{g,k}$. 

(ii$'$) The image of $h$ consists of a single point of $\Cal M_{g,k}$ 
(i.e., $h: S\to\Cal M_{g,k}$ is set-theoretically constant). 

(iii) For each generic point $\eta$ of $S$, $h(\eta)$ is a closed point. 

(iv) For each $1$-dimensional irreducible, reduced, closed subscheme $C$ of $S$, 
the 
%(relative) 
$C$-curve $f_C: X\times_SC\to C$ is isotrivial. 

(iv$'$) For each irreducible component $W$ of $S$ (regarded as a reduced closed subscheme of $S$), 
there exists a closed point $s$ of $W$, such that, 
for each $1$-dimensional irreducible, reduced, closed subscheme $C$ of $W$ 
passing through $s$, 
the 
%(relative) 
$C$-curve $f_C: X\times_SC\to C$ is isotrivial. 
\endproclaim

\demo{Proof}
Standard. 
\qed\enddemo

\proclaim{Theorem 3.9} Assume that 
%$g\geq 2$ and that 
$f:X\to S$ is non-isotrivial. 
Then $\pi_1$ is not constant on $S$. 
\endproclaim 

In the process of proving Theorem 3.9 we prove the following more precise result.

\proclaim{Theorem 3.10} 
Assume that 
%$g\geq 2$ and that 
$f:X\to S$ is non-isotrivial. 
Then there 
exist a connected finite \'etale cover $S'\to S$, 
a connected finite \'etale cover $Y'\to X'\defeq X\times _SS'$ 
whose Galois closure is geometrically connected over $S'$ and 
of degree prime to $p$ over $X\times_SS'$, 
a generic point $\eta'\in S'$ and 
a %closed 
point $v'\in \overline{\{\eta'\}}\subset S'$, 
such that the $p$-rank of the geometric fibre 
$Y'_{\bar v'}$ of $Y'$ above $v'$ is smaller than 
the $p$-rank of the geometric fibre 
$Y'_{\bar \eta'}$ of $Y'$ above $\eta'$ 
(thus, $v'\neq\eta'$ necessarily). 
\endproclaim

%Theorem 3.9 implies the following where we consider the full tame fundamental group 
%(cf. [S]. Theorem 4.8, and the proof therein that Theorem 3.9 implies Theorem 3.12).
Theorem 3.9 admits certain variants 
(Theorems 3.12 and 3.13) for curves which are not necessarily proper. 
To state them, let 
$X$ be a proper and smooth $S$-curve of (constant) genus $g\geq 0$ and 
$D\subset X$ a relatively \'etale divisor of (constant) degree $n\geq 0$. 
Given a point $s\in S$, choose 
a geometric point $\bar s$ above $s$ and let $(X_{\bar s}, D_{\bar s})$ be the 
geometric fibre of $(X,D)/S$ at $\bar s$. 
Then the isomorphism type of 
the (geometric) tame fundamental group $(\pi_1^t)_s\defeq \pi_1^t(X_{\bar s}-D_{\bar s})$ 
(resp. the (geometric) fundamental group $(\pi_1)_s\defeq \pi_1(X_{\bar s}-D_{\bar s})$) 
is independent of (resp. (if $n>0$) dependent on) the choice of $\bar s$.  
%(and the implicit base point on $X_{\bar s}$ used to define $\pi_1(X_{\bar s})$). 
%If $s$ and $t$ are points in $S$ such that $s\in \overline {\{t\}}$ holds, 
%then Grothendieck's specialisation theory for fundamental groups
%implies the existence of a (continuous) surjective specialisation homomorphism 
%%\Sp=
%\Sp_{t,s}:\pi_1(X_{\bar t})\twoheadrightarrow \pi_1(X_{\bar s})$. 

We say that $(X,D)/S$ is {\bf isotrivial}, if there exist 
a finite extension $k'/k$, a connected $k'$-scheme $S'$, a finite \'etale $k'$-morphism 
$S'\to S\times_kk'$, a proper and smooth $k'$-curve $X_0'$ and a relatively \'etale divisor 
$D_0'$ of $X_0'$, such that $(X,D)\times_S S'$ is isomorphic to $(X_0',D_0')\times_{k'}S'$ 
over $S'$. 

\demo{Remark 3.11}
(i) As in the proper case, 
if $s$ and $t$ are points in $S$ such that $s\in \overline {\{t\}}$ holds, 
then Grothendieck's specialisation theory for tame fundamental groups
implies the existence of a (continuous) surjective specialisation homomorphism 
$%\Sp=
\Sp^t_{t,s}:\pi_1^t(X_{\bar t}-D_{\bar t})\twoheadrightarrow 
\pi_1^t(X_{\bar s}-D_{\bar s})$. 
Thus, by Definition 3.1, we have the notion of $\Sp^t$-constancy of $(X,D)/S$ 
as well, which is equivalent 
to the notion of constancy of $(X,D)/S$ by Lemma 3.2. (Note that the tame 
fundamental group of affine curves is finitely generated (cf. [SGA1], Expos\'e XIII, Corollaire 2.12).) 
Note that no such specialisation homomorphisms are available for (full) 
fundamental groups, if $n>0$. 

(ii) Assume that $2-2g-n<0$. Then, 
as in the proper case (cf. Lemma 3.8), 
$(X,D)/S$ is isotrivial if and only if 
the (coarse) classifying morphism $h: S\to\Cal M_{g,[n],k}$ for $(X,D)/S$ 
is set-theoretically constant. 
Here, $\Cal M_{g,[n],k}\ (=\Cal M_{g,[n],\Bbb F_p}\times _{\Bbb F_p}k)$ is the 
coarse moduli space of proper, smooth and geometrically connected curves of genus $g$ 
equipped with a relatively \'etale divisor of degree $n$ over $k$-schemes. 
\enddemo

\proclaim{Theorem 3.12} 
%Let $k$ be a field of characteristic $p$. Let $S$ be a smooth and irreducible affine $k$-curve,
%and $f:X\to S$ a {\bf non-isotrivial} proper and smooth 
%%(relative) 
%$S$-curve of (constant) genus $g$. 
%Let $\{s_1,s_2,...,s_n\}$ be $n$-sections of $f$ with disjoint support such that $2-2g-n<0$. 
%Then the tame fundamental group is not constant on the pair $(f,\{s_1,s_2,...,s_n\})$, i.e., there
%exists a closed point $s\in S$, $\eta \in S$ is the generic point,  
%such that the specialisation homomorphism 
%$\Sp:\pi_1^t(X_{\bar \eta}-\{s_1(\bar \eta),s_2(\bar \eta),...,s_n(\bar \eta)\})
%\to \pi_1^t(X_{\bar s}-\{s_1(\bar s),s_2(\bar s),...,s_n(\bar s)\})$ between tame fundamental groups, 
%where $X_{\bar \eta}\defeq X\times _S\bar \eta$ and $X_{\bar s}\defeq X\times _S\bar s$, 
%is not an isomorphism.
Assume that $2-2g-n<0$ and that $(X,D)/S$ is non-isotrivial. 
Then $\pi_1^t$ is not constant on $S$. 
\endproclaim

\proclaim{Theorem 3.13} 
Assume that $2-2g-n<0$ and that $(X,D)/S$ is non-isotrivial. 
Then $\pi_1$ is not constant on $S$. 
\endproclaim

\subhead
\S 4. Proof of the Main Theorems
\endsubhead
In this section we prove the main results: Theorems 
3.6, 3.9, 3.10, 3.12 and 3.13. 
First, we work with the assumptions in Theorems 3.9 and 3.10. 
In particular, 
$k$ is a field of characteristic $p>0$,
$S$ is a connected and reduced $k$-scheme of finite type, 
and $f:X\to S$ is a non-isotrivial proper and smooth 
%(relative) 
$S$-curve of genus $g\geq 2$. 

\demo{Proof of Theorem 3.10} 

\proclaim{Lemma 4.1} 
Let $s$ be a point of $S$ and $X_{\bar s}$ the geometric fibre of $f:X\to S$ 
at a geometric point $\bar s$ above $s$. Let $Y_{\bar s}\to X_{\bar s}$ be 
a finite \'etale cover whose Galois closure is of degree prime to $p$ over 
$X_{\bar s}$. Then there exist a connected finite \'etale cover $S'\to S$ and 
a connected finite \'etale cover $Y'\to X'\defeq X\times _SS'$ 
whose Galois closure is geometrically connected over $S'$ and 
of degree prime to $p$ over $X\times_SS'$, 
a point $s'\in S'$ above $s\in S$ and a geometric point $\bar s'$ above $s'$ 
which dominates $\bar s$, 
such that the cover $(Y')_{\bar s'} \to (X')_{\bar s'}$ is isomorphic to 
the pull-back of $Y_{\bar s}\to X_{\bar s}$ to $\bar s'$. 
\endproclaim

\demo{Proof} 
%(See also [S], Lemma/Definition 4.7, when $Y_{\bar s}\to X_{\bar s}$ is a 
%$\mu_N$-torsor for some positive integer $N$ prime to $p$.) 
This follows from the fact (cf. [St], Proposition 2.3) that the natural sequence of 
profinite groups 
$$
1
\to \pi_1(X_{\bar s})^{p'}
\to \pi_1(X)^{(p')}
\to \pi_1(S) 
\to 1
$$
%arising from the homomorphisms of profinite groups 
%$\pi_1(X_{\bar s})\to \pi_1(X)$ and 
%$\pi_1(X)\to \pi_1(S)$ 
%associated (by functoriality of $\pi_1$) to 
%the natural morphisms of schemes 
%$X_{\bar s}\to X$ and 
%$X\to S$, respectively, 
is exact, where $\pi_1(X_{\bar s})^{p'}$ stands for the maximal pro-prime-to-$p$ quotient 
of $\pi_1(X_{\bar s})$ and $\pi_1(X)^{(p')}$ stands for the maximal quotient of 
$\pi_1(X)$ in which the image of $\Ker(\pi_1(X)\to \pi_1(S))$ is pro-prime-to-$p$. 
\qed\enddemo

%By Lemma 4.1, we may reduce the problem to the case where $S$ is 
%a smooth, separated, geometrically connected curve over $k$. 
%Indeed, by Lemma 3.8, there exists a $1$-dimensional, irreducible, reduced, 
%closed subscheme $C$ of $S$, such that $f_C: X\times_SC\to C$ is non-isotrivial. 
%Next, there exist a finite extension $k_1$ of $k$, 
%a smooth, geometrically connected curve $C_1$ over $k_1$ and 
%a dominant $k$-morphism $C_1\to C$. (For example, take 
%a connected component of the regular locus of $(C\times \bar k)_{\text{red}}$ and 
%consider a descent to a finite extension of $k$.) As $C_1\to C$ is dominant, 
%$f_{C_1}: X\times_SC_1\to C_1$ is non-isotrivial. Now, suppose that the assertion 
%of Theorem 3.10 holds with $X\overset{f}\to{\to} S\to \Spec(k)$ replaced by 
%$X\times_SC_1\overset{f_{C_1}}\to{\to} C_1\to\Spec(k_1)$. 
%Then there exist a connected finite \'etale cover $C_1'\to C_1$, 
%a connected finite \'etale cover $Y'\to X\times _SC_1'$ 
%whose Galois closure is geometrically connected over $C_1'$ and 
%of degree prime to $p$ over $X\times_SC_1'$, 
%a closed point $c'\in C_1'$ and a generic point $\xi'\in C_1'$, 
%such that the $p$-rank of the geometric fibre 
%$Y'_{\bar c'}$ of $Y'$ above $c'$ is smaller than 
%the $p$-rank of the geometric fibre 
%$Y'_{\bar \xi'}$ of $Y'$ above $\xi'$. 

We may reduce the problem to the case where $S$ is geometrically connected and 
geometrically reduced over $k$. 
% and each irreducible component of $S$ is geometrically irreducible over $k$. 
Indeed, first, take 
a connected component of $S\times_kk^{\text{sep}}$ 
(regarded 
%as an open and closed subscheme and 
as a scheme over $k^{\text{sep}}$), 
which descends to a scheme $S_1$ over a finite separable extension $k_1$ 
of $k$. 
Next, 
%consider a finite separable extension $k_2$ of $k_1$ over which 
%each irreducible component of $S_1\times_{k_1}k_1^{\text{sep}}$ is defined 
%and set $S_2\defeq S_1\times_{k_1}k_2$. 
%Finally, 
consider the reduced closed subscheme 
$(S_1\times_{k_1}k_1^{\text{perf}})_{\text{red}}$ of 
$S_1\times_{k_1}k_1^{\text{perf}}$ 
(regarded as a scheme over $k_1^{\text{perf}}$), which  
descends to a scheme $S_2$ over a 
finite purely inseparable extension $k_2$ of $k_1$. 
Then $S_1$ is geometrically connected over $k_1$, 
%$S_2$ is geometrically connected over $k_2$ and 
%each irreducible component of $S_2$ is geometrically irreducible, 
and 
$S_2$ is geometrically connected and geometrically reduced over $k_2$. 
As $S_1$ is a connected finite \'etale cover of $S$, we may replace 
$X\to S\to \Spec(k)$ by $X\times_SS_1\to S_1\to \Spec(k_1)$. 
As the morphisms $S_2\to S_1$ and $X\times_SS_2\to X\times_SS_1$ 
preserve the categories of finite \'etale covers, we may replace 
$X\times_SS_1\to S_1\to \Spec(k_1)$ by $X\times_SS_2\to S_2\to \Spec(k_2)$. 

So, from now on, we will assume that $S$ 
is geometrically connected and geometrically reduced over $k$.
% and that each irreducible component of $S$ is geometrically irreducible. 

\proclaim {Lemma 4.2} 
%(i)\ The morphism $f:X\to S$ descends to a finitely generated field $k_0$. 
%More precisely, there exists a finitely generated 
%(over the prime field $\Bbb F_p$) field $k_0\subset k $, an affine, smooth, and irreducible 
%$k_0$-curve $S_0$, a proper and smooth morphism $f_0:X_0\to S_0$ such that 
%$S=S_0\times _{k_0}k$, $X=X_0\times _{k_0}k$, and we have a commutative diagram with cartesian squares 
%$$
%\CD
%X @>{f}>> S @>>> \Spec (k) \\
%@VVV    @VVV  @VVV \\
%X_0  @>{f_0}>>  S_0  @>>> \xi\defeq \Spec (k_0) \\
%\endCD
%$$
%where the right vertical map is the natural morphism. 
%
%(ii)\ There exists a finitely generated $\Bbb F_p$-algebra $R$, with $\Fr (R)=k_0$, a smooth 
%equidimensional morphism $\Cal S\to T\defeq \Spec (R)$ of relative dimension $1$
%with $\Cal S$ connected, 
%a proper and smooth 
%%(relative) 
%$\Cal S$-curve $\tilde f:\Cal X\to \Cal S$ of (constant) 
%genus $g$, such that we have a commutative diagram with cartesian squares
%$$
%\CD
%\Cal X @>{\tilde f}>> \Cal S @>>> T=\Spec (R) \\
%@AAA     @AAA   @AAA \\
%X_0  @>{f_0}>>  S_0  @>>> \xi\defeq \Spec (k_0) \\
%\endCD
%$$
%where the right vertical map is the natural morphism induced by the inclusion $R\subset k_0$. 
%%Thus, $f_0$ is the generic fibre of $\tilde f$. 
There exist a finitely generated $\Bbb F_p$-subalgebra $R$ of $k$, 
a connected scheme $\Cal S$ flat, of finite type
over $R$ with geometrically connected and geometrically reduced fibres,  
%such that each irreducible component of any fibre of $\Cal S\to\Spec(R)$ 
%is geometrically irreducible, 
and a proper and smooth 
%(relative) 
$\Cal S$-curve $\tilde f:\Cal X\to \Cal S$ of (constant) 
genus $g$, such that we have a commutative diagram with cartesian squares
$$
\CD
\Cal X @>{\tilde f}>> \Cal S @>>> \Spec (R) \\
@AAA     @AAA   @AAA \\
X  @>f>>  S  @>>> \Spec (k) \\
\endCD
\tag{$\ast$}
$$
where the right vertical map is the natural morphism induced by the inclusion $R\subset k$. 
\endproclaim

\demo{Proof} 
As $k$ is the direct limit of finitely generated $\Bbb F_p$-subalgebras, 
it follows from [EGA IV], Th\'eor\`eme (8.8.2) that 
there exists a finitely generated $\Bbb F_p$-subalgebra $R$ of $k$, 
%a reduced scheme $\Cal S$ flat, of finite type over $R$ with connected fibres, 
schemes $\Cal S$ and $\Cal X$ of finite type over $R$, and 
an $R$-morphism $\tilde f:\Cal X\to \Cal S$,
such that we have a commutative diagram ($\ast$) with cartesian squares. 
Replacing $R$ by a finitely generated $\Bbb F_p$-subalgebra $R'$ of $k$ containing $R$ 
and $\tilde f:\Cal X\to \Cal S$ by $\tilde f_{R'}=\tilde f\times_RR': 
\Cal X\times_RR'\to \Cal S\times_RR'$, we may assume that 
$\Cal S$ is flat over $R$ ([EGA IV], Th\'eor\`eme (11.2.6)) 
with geometrically connected and geometrically reduced fibres 
([EGA IV], Th\'eor\`eme (9.7.7)) 
%and  
%each irreducible component of any fibre of $\Cal S\to\Spec(R)$ 
%is geometrically irreducible ([EGA IV], Lemme (9.7.1) and Proposition (9.7.8)) 
and that $\Cal X$ is proper and smooth 
over $\Cal S$ ([EGA IV], Th\'eor\`eme (8.10.5) and Proposition (17.7.8)). 
As $\Cal S\to \Spec(R)$ is flat of finite type, it is (universally) open 
([EGA IV], Th\'eor\`eme (2.4.6)), hence its generic fiber, which is (geometrically) 
connected, is dense. This implies that $\Cal S$ is connected. 
As $\Cal X\to\Cal S$ is proper and smooth and has 
geometrically connected fibres on the image of $S$ in $\Cal S$, 
we conclude (by observing the Stein factorisation) that 
$\Cal X\to \Cal S$ must have geometrically connected fibres everywhere 
(cf. [EGA III], Remarque (7.8.10) and [SGA1], Expos\'e X, Proposition 1.2). 
Finally, as $\Cal X\to \Cal S$ is proper and flat and $\Cal S$ is connected, 
the dimension and the (arithmetic) genus of the fibres are constant. 
Thus, $\Cal X \to \Cal S$ is a proper and smooth 
%(relative) 
$\Cal S$-curve 
of constant genus $g$, as desired. 
\qed
\enddemo

Next, set $T\defeq\Spec(R)$ 
and 
let $\xi$ 
%and $t$ 
be the generic point 
%and a closed point 
of $T$. 
%, respectively. 
%Thus, the residue field $k(t)$ at $t$ is a finite field of characteristic $p$. 
For each $t\in T$, set $\Cal S_t\defeq \Cal S\times _T t$ and 
$\Cal X_t\defeq \Cal X\times _Tt$. 
We have a commutative diagram with cartesian squares 

$$
\CD
\Cal X @>{\tilde f}>> \Cal S @>>> T=\Spec (R) \\
@AAA     @AAA   @AAA\\
\Cal X_t  @>{\tilde f_t}>>  \Cal S_t  @>>> t=\Spec (\kappa(t)) \\
\endCD
$$
where $\Cal S_t\to \Spec (\kappa(t))$ is of finite type, 
geometrically connected and geometrically reduced, 
and $\tilde f_t:\Cal X_t\to \Cal S_t$ is a proper and smooth 
%(relative) 
$\Cal S_t$-curve of (constant) genus $g$. 

\proclaim {Lemma 4.3} 
There exists a closed point $t\in T$ such that 
the 
%(relative) 
$\Cal S_t$-curve $\tilde f_t:\Cal X_t\to \Cal S_t$ is non-isotrivial.
\endproclaim

\demo{Proof} 
Let 
$h:S\to\Cal M_{g,k}$ 
and 
$\tilde h:\Cal S \to\Cal M_{g,\Bbb F_p}$ 
be the (coarse) classifying morphisms for 
$f: X\to S\ (\to\Spec(k))$ and
$\tilde f: \Cal X\to \Cal S\ (\to\Spec(\Bbb F_p))$. 
Let 
$\tilde\beta: \Cal S\to T$ 
be the structure morphism and 
set $\tilde h_T \defeq (\tilde h, \tilde\beta): \Cal S \to \Cal M_{g,\Bbb F_p}\times_{\Bbb F_p} T$. 
For each $t\in T$, let $\tilde h_t: \Cal S_t\to \Cal M_{g,\kappa(t)}$ be the 
(coarse) classifying morphism for 
$\tilde f_t: \Cal X_t\to \Cal S_t\  (\to\Spec(\kappa(t)))$. 
Then we have $\tilde h_t=\tilde h_T \times_Tt$. Further, 
we have $\tilde h_T \times_Tk=(\tilde h_\xi)\times_{\kappa(\xi)}k=h$. 

By Lemma 3.8, there exists a generic point $\eta$ of $S$ 
such that $h(\eta)$ is not a closed point of $\Cal M_{g,k}$. 
Let $\tilde \eta$ be the image of $\eta$ in $\Cal S$, which is a generic point of $\Cal S$. 
Then 
%$\xi\defeq \tilde\beta(\tilde\eta)$ is the generic point of $T$ 
$\tilde\beta(\tilde\eta)=\xi$, 
and 
$\tilde h_{\xi}(\tilde\eta)$ is not a closed point of $\Cal M_{g,\kappa(\xi)}$. 
Note that there exists an open neighbourhood $\Cal W$ of $\tilde \eta$ which is irreducible. 
(For example, remove all the irreducible components of $\Cal S$ but $\overline{\{\tilde\eta\}}$.) 
Set $\alpha\defeq \tilde h_T |_{\Cal W}$ and $\beta\defeq \tilde\beta|_{\Cal W}$. 
%It suffices to prove that there exists a closed point $t\in T$ such that 
%$\tilde h|_{\Cal W_t}$ is non-constant. 
By [EGA IV], Corollaire (9.2.6.2), there exists 
a non-empty open subset $\Cal U$ of 
$\Cal W$ such that $\dim(\alpha^{-1}(\alpha(s)))$ and 
$\dim(\beta^{-1}(\beta(s)))$ are constant for $s\in \Cal U$. 
Finally, for each $t\in T$, set $\Cal W_t\defeq\Cal W\times_Tt$. 

Now, take any closed point $s$ of $\Cal U$ and set $t\defeq \beta(s)$, 
which is a closed point of $T$. 
As $\tilde h_{\xi}|_{\Cal W_{\xi}}$ 
(or, equivalently, $\tilde h|_{\Cal W_{\xi}}$) 
is non-constant and $\Cal W_{\xi}$ is irreducible, we have 
$$
\dim(\alpha^{-1}(\alpha(\tilde\eta)))
<\dim(\beta^{-1}(\beta(\tilde\eta)))
=\dim(\Cal W_{\xi}),$$
hence 
$$
\dim(\alpha^{-1}(\alpha(s))
<\dim(\beta^{-1}(\beta(s)))
=\dim(\Cal W_t), 
$$
which implies that 
%$\tilde h|_{\Cal W_t}$ (or, equivalently, $\tilde h_t|_{\Cal W_t}$) is non-constant, 
$\tilde h$ (or, equivalently, $\tilde h_t$) is non-constant on each irreducible component 
of $\Cal W_t$ of maximal dimension ($=\dim(\Cal W_t)$). 
As $\Cal M_{g,\kappa(t)}\to \Cal M_{g,\Bbb F_p}$ is finite, 
this implies that 
$\tilde h_t|_{\Cal S_t}$ is non-constant, 
as desired. 
\qed
\enddemo

Next, consider the (relative) curve $\Cal X_1\to \Cal S$ 
where $\Cal X_1$ is the Frobenius twist of $\Cal X$ (cf. $\S1$), 
$\Cal J_1\defeq \Pic^0(\Cal X_1/\Cal S)$ the (relative) Jacobian of $\Cal X_1$ 
over $\Cal S$ 
which is an $\Cal S$-abelian scheme,
and $\Theta \hookrightarrow \Cal J_1$ the (relative) Raynaud theta divisor 
for $\Cal X\to\Cal S$ 
(cf. $\S1$). 
Let $t$ be a closed point of $T$ such that the relative curve $\Cal X_t\to \Cal S_t$ is 
%{\bf 
non-isotrivial 
%} 
(cf. Lemma 4.3). 
Set $\Cal J_{1,t}\defeq \Cal J_1\times _Tt$ and $\Theta _t\defeq \Theta \times _Tt$, 
which are the (relative) Jacobian 
of $\Cal X_{1,t}\defeq \Cal X_1\times _Tt$ over $\Cal S_t$ 
and the (relative) Raynaud theta divisor 
for 
%the relative curve 
$\Cal X_t\to \Cal S_t$, respectively. 

Let $\Cal S^{\text{sm}}\subset \Cal S$ and $\Cal S_t^{\text{sm}}\subset \Cal S_t$ 
be the smooth loci for $\Cal S \to T$ and $\Cal S_t \to t$, respectively. 
As $\Cal S \to T$ is of finite type, 
$\Cal S^{\text{sm}}\subset \Cal S$ and $\Cal S_t^{\text{sm}}\subset \Cal S_t$ are open. 
As $\Cal S \to T$ is flat, 
we have	 $(\Cal S^{\text{sm}})_t=\Cal S_t^{\text{sm}}$. 
As $\Cal S \to T$ has geometrically reduced fibres, 
$\Cal S^{\text{sm}}\subset \Cal S$ is dense in each fibre, and, 
in particular, $\Cal S_t^{\text{sm}}$ is dense in $\Cal S_t$. 
Let $\Cal S_t^0$ be the (disjoint) union of connected 
(or, equivalently, irreducible) components 
$V$ of $\Cal S_t^{\text{sm}}\subset \Cal S_t$ for which 
the $V$-curve $\tilde f_{V}= (\tilde f_t)_{V}: 
\Cal X\times_TV=\Cal X_t\times_tV\to V$ is non-isotrivial. 
By Lemma 4.3, $\Cal S_t^0$ is non-empty. 

Take any component $V$ of $\Cal S_t^0$ and 
any closed point $s$ of $V$. Then, 
by Lemma 3.8, 
there exists a $1$-dimensional, irreducible, closed 
subscheme $C$ of $V$ passing through $s$ 
such that the 
%(relative) 
$C$-curve 
$\tilde f_C=(\tilde f_V)_C: \Cal X\times_TC\to C$ is non-isotrivial. 
Let $\gamma$ and $u$ be the generic points of $C$ and $V$, 
respectively. 
Then it follows from Theorem 2.1 that the specialisation homomorphism 
$\Sp:\pi_1(\Cal X_{t,\bar \gamma})\to \pi_1(\Cal X_{t,\bar s})$ is not an isomorphism. 
Here, $\bar \gamma$ (resp. $\bar s$) is a geometric point above $\gamma$ (resp. $s$), and 
$\Cal X_{t,\bar \gamma}\defeq \Cal X_t\times _t \bar \gamma$ 
(resp. $\Cal X_{t,\bar s}\defeq \Cal X_t\times _s \bar s$). 

%We recall the following technical result (cf. [S], Lemma/Definition 4.7).
%\proclaim{Proposition 3.8} Let $S$ be an $\Bbb F_p$-scheme, and $X\to S$ a proper and smooth 
%%(relative) 
%$S$-curve. Then there exists a positive integer $d$, such that the following holds: 
%Let $N$ be a positive integer prime to $pd$. 
%Let $s$ be a closed point of $S$, and $f_s:Y_s\to X_s\defeq X\times _Sk(s)$ 
%%an \'etale 
%a $\mu_N$-torsor.  
%%over the fibre of $X$ above the point $s$. 
%%If $N$ is coprime to $d$, then 
%Then there exists a finite \'etale cover $h:S'\to S$, 
%an %\'etale 
%$\mu_N$-torsor $f':Y'\to X'\defeq X\times _SS'$, 
%and a closed point $s'$ of $S'$ with $h(s')=s$, 
%such that the geometric fibre $f'_{\bar s'}:Y'_{\bar s'}\to X'_{\bar s'}$ of 
%the torsor $f'$ above the point $s'$ is isomorphic (as a $\mu_N$-torsor) 
%to the torsor $f_{\bar s}\defeq f_s\times _{s}{\bar s}:Y_{\bar s}\to X_{\bar s}$.
%\endproclaim

%More precisely, 
Now, 
it follows from Proposition 2.4 and Lemma 4.1 
that the following hold. There exist a finite \'etale cover $\Cal S'\to \Cal S$, 
a finite \'etale cover $\Cal Y'\to \Cal X'\defeq \Cal X\times _{\Cal S}\Cal S'$ 
whose Galois closure is 
of degree prime to $p$, 
a (generic) point $u'\in \Cal S'_t\defeq \Cal S'\times _Tt$ above $u$, 
a point $\gamma'\in\overline{\{u'\}}\subset\Cal S'_t$ above $\gamma$ 
and 
a (closed) point $s'\in \overline{\{\gamma'\}}\subset{\Cal S'_t}$ above $s$,  
such that
$\sat (\Theta '_{\bar u'}\{p'\})
\subset
\sat (\Theta '_{\bar \gamma'}\{p'\})
\subsetneq \sat (\Theta' _{\bar s'}\{p'\})$ holds 
in $\Cal J'_{1}\{p'\}\defeq \Cal J'_{1,\bar u'}\{p'\}
=\Cal J'_{1,\bar \gamma'}\{p'\}
=\Cal J'_{1,\bar s'}\{p'\}$, where 
$\Cal J'_1\defeq \Pic^0(\Cal Y'_1/\Cal S')$ is the (relative) Jacobian of $\Cal Y'_1$ over $\Cal S'$, 
$\Theta'\hookrightarrow \Cal J'_1$ is the (relative) Raynaud theta divisor for 
$\Cal Y'\to\Cal S'$, 
and $\bar u'$ (resp. $\bar\gamma'$, resp. $\bar s'$) 
is a geometric point above $u'$ (resp. $\gamma'$, resp. $s'$).
In particular, there exists $x\in \Cal J'_{1}\{p'\}$ 
such that $\sat (x)\cap \Theta'_{\bar u'}\{p'\}=\varnothing$ and $\sat (x)\cap \Theta'_{\bar s'}\{p'\}
\neq \varnothing$.  Let $z\in \sat (x)\cap \Theta'_{\bar s'}\{p'\}$, and $N\defeq \ord(z)$. Let $\Cal J_1'[N]\defeq \Ker (\Cal J_1'@>{[N]}>> \Cal J_1')$ be the kernel of multiplication by $N$ on the abelian scheme $\Cal J_1'$. Thus,
$\Cal J_1'[N]$ is a finite \'etale commutative $\Cal S'$-group scheme which is \'etale-locally isomorphic to 
$(\Bbb Z/N\Bbb Z)^{2g}$. 
After possibly passing to a finite \'etale cover 
%$\Cal S''\to \Cal S'$ 
of $\Cal S'$ 
which 
``trivialises'' $\Cal J_1'[N]$, we can assume (without loss of generality) that there exists a section 
$\sigma:\Cal S'\to \Cal J_1'[N]$ of the natural projection $\Cal J_1'[N]\twoheadrightarrow \Cal S'$, with image $\Cal Z\defeq \Im (\sigma)$, such that
$\Cal Z_{\bar s'}\defeq \Cal Z\times _{\Cal S'}\bar s'=z$ holds. Write 
$\Cal Z.\Theta'$ 
for the scheme-theoretic  intersection of $\Cal Z$ and $\Theta '$ inside $\Cal J_1'$. 
By definition, we have $\sigma(s')\in \Cal Z.\Theta'\neq\varnothing$. 
(By a slight abuse of notation, we write $z=\sigma(s')$.) 
We have natural morphisms 
$\Cal Z.\Theta'\to \Cal S'\to \Cal S\to T$. 
The following is the main technical ingredient of our proof. 

\proclaim {Proposition 4.4} The morphism 
$\Cal Z.\Theta'\to T$ is flat at $z\in \Cal Z.\Theta'$. 
\endproclaim

\demo{Proof} 
%For a closed point $s'\in \Cal S'$, write $s\in \Cal S$ (resp. $t\in T$) 
%for its image in $\Cal S$ (resp. in $T$), 
%$z\in \Cal Z$ for its pre-image in $\Cal Z$ via the isomorphism $\Cal Z\isom \Cal S'$, 
%$\Cal O_{\Cal S',s'}$, $\Cal O_{\Cal S,s}$, $\Cal O_{T,t}$, and $\Cal O_{\Cal Z,z}$ 
%for the corresponding local rings. 
%Let $z$ be any closed point of $\Cal Z.\Theta'\subset \Cal Z$ above $t$, and 
%$s'$ and $s$ the images of $z$ in $\Cal S'$ and $\Cal S$, respectively. 
We have natural morphisms 
$\Cal O_{T,t}\to \Cal O_{\Cal S,s}\to \Cal O_{\Cal S',s'}=\Cal O_{\Cal Z,z}$. 
Note that $\Cal O_{\Cal S',s'}$ is $\Cal O_{T,t}$-flat, since 
$\Cal S'$ is \'etale over $\Cal S$ and $\Cal S$ is 
%smooth 
flat 
over $T$. Moreover, there exists a natural surjective homomorphism
$\Cal O_{\Cal J_1',z}\twoheadrightarrow \Cal O_{\Cal Z,z}$. Let 
$\Cal I_{\Theta',z}=(f)\subset \Cal O_{\Cal J_1',z}$ be the ideal defining the theta divisor $\Theta'$ 
(cf. Theorem 1.1). We will show that $\Cal O_{\Cal Z.\Theta',z}=
{\Cal O_{\Cal Z,z}}/(f)$ is flat over $\Cal O_{T,t}$.

Let $\kappa(t)$ be the residue field of $T$ at $t$, and $g$ the image of $f$ 
in $M\defeq \Cal O_{\Cal Z,z}\otimes _{\Cal O_{T,t}}\kappa(t)$.
By [EGA IV], Chapitre 0, Proposition (15.1.16), $c)\implies a)$, 
it suffices to show  that $g$ is $M$-regular. 
Furthermore, by [EGA IV], Chapitre 0, Proposition (16.3.6) and (Proposition) 16.5.5, $b)\implies a)$, 
to show the latter it suffices to show that $\dim (M/gM)<\dim M$.  
Observe that $M=\Cal O_{\Cal Z,z}\otimes _{\Cal O_{T,t}}\kappa(t)=\Cal O_{\Cal Z_t,z}$ 
and that $M/gM= \Cal O_{(\Cal Z.\Theta)_t,z}$, 
where $\Cal Z_t\defeq \Cal Z\times_Tt$. 
%and $\Cal S'_t\defeq \Cal S'\times _{T}t$. 
%is a discrete valuation ring since $\Cal S'_t\to t=\Spec (\kappa(t))$ is a proper and smooth $\kappa(t)$-curve. 
%Furthermore, $\dim (M/gM)=\dim \Cal O_{(\Cal Z.\Theta)_t,z}=0$ occurs by our assumption and 
%choice of $\Cal Z$ that $(\Cal Z.\Theta)_{t,\bar u'}=\varnothing$ (cf. discussion before Proposition 4.4).
Now, since $(\Cal Z.\Theta)_{t,\bar u'}=\varnothing$ (cf. discussion before Proposition 4.4), 
we have $\dim (M/gM)<\dim M$, as desired. 
\qed
\enddemo

%Let $S_0'\defeq S_0\times _{\Cal S}\Cal S'$ and $X_0'\defeq X_0\times _{S_0}S_0'$. 
Set 
$S_0\defeq \Cal S\times_T\xi$,  
$S_0'\defeq \Cal S'\times_T\xi$, 
$X_0\defeq \Cal X\times_T\xi$, 
and $X_0'\defeq X_0\times _{S_0}S_0'=(\Cal X\times_{\Cal S}\Cal S')\times_T\xi$. 
%We will show that the (geometric) fundamental group $\pi_1$ is not constant 
%on the relative curve $X_0'\to S_0'$, which easily implies 
%that $\pi_1$ is not constant on the relative curve
%$X_0\to S_0$. 
%In the next discussion we will write $z\in \Cal J_1'\{p'\}$ 
%for the element corresponding to the above section $\Cal Z\hookrightarrow \Cal J_1'[N]$.
Let $\eta'$ be the generic point of $S_0'$. Then $\sat (z) \cap \Theta'_{\bar \eta'}=\varnothing$ since $\sat (z) \cap \Theta'_{\bar u'}   =\varnothing$ and $\eta'$ specialises 
into $u'$. Next, we claim that there exists a 
point $v'\in {S_0'}$ such that $\sat (z) \cap \Theta'_{\bar v'}\neq \varnothing$, 
which would imply 
%that the specialisation homomorphism  $\pi_1((X_0')_{\bar \eta'})\to \pi_1((X_0')_{\bar v'})$ 
%is not an isomorphism (cf. Lemma 2.3). 
the assertion of Theorem 3.10 (cf. Proposition 1.3). 
%Recall $\xi=\Spec (k_0)$ is the generic point of $T$ (cf. Notations in Lemma 4.2). 
As the morphism $\Cal Z.\Theta'\to T$ is flat at $z\in \Cal Z.\Theta'$ 
(which is above $t\in T$), 
its image contains the image of the natural morphism $\Spec(\Cal O_{T,t})\to T$, 
hence, in particular, contains the generic point $\xi\in T$. 
So, let $y\in \Cal Z.\Theta'$ be any point above $\xi\in T$, and $v'\in \Cal S'$ 
the image of $y$ in $\Cal S'$. Then 
%$v'$ is a closed point of $S_0'$, since $v'$ maps to $\xi$ in $T$, and 
by assumption it holds that  $z\in \sat (z) \cap \Theta'_{\bar v'}\neq \varnothing$ as claimed. 

This finishes the proof of Theorem 3.10. 
\qed
\enddemo

\demo{Proof of Theorem 3.9} 
%Let the notations as in the proof of Theorem 3.10. 
%Let $\eta$ be the generic point of $S_0$, and $v$ the image of $v'$ in $S_0$, 
%then it follows from Theorem 3.10 and Lemma 2.3 that the specialisation homomorphism 
%$\pi_1((X_0)_{\bar \eta})\to \pi_1((X_0)_{\bar v})$ is not an isomorphism, 
%which implies that the (geometric) fundamental group $\pi_1$ is not constant 
%on the (relative) curve $X_0\to S_0$. Hence the (geometric) fundamental group $\pi_1$ 
%is not constant on the (relative) curve $X\to S$.
%
%This finishes the proof of Theorem 3.9. 
Theorem 3.9 follows easily from Theorem 3.10. Indeed, notations as in Theorem 3.10, 
write $\eta$ and $v$ for the images of $\eta'$ and $v'$ in $\Cal S$, respectively. 
Then we have the following commutative diagram 
$$
\matrix
\pi_1(Y'_{\bar\eta'})&\twoheadrightarrow&\pi_1(Y'_{\bar v'})\\
&&\\
\cap&&\cap\\
&&\\
\pi_1(X'_{\bar\eta'})&\twoheadrightarrow&\pi_1(X'_{\bar v'})\\
&&\\
\Vert&&\Vert\\
&&\\
\pi_1(X_{\bar\eta})&\twoheadrightarrow&\pi_1(X_{\bar v}), 
\endmatrix
$$
where the horizontal arrows are specialisation homomorphisms. The assertion of Theorem 3.10 
implies that $\Sp: \pi_1(Y'_{\bar\eta'})\twoheadrightarrow \pi_1(Y'_{\bar v'})$ is not  
injective, which implies that 
$\Sp: \pi_1(X_{\bar\eta})\twoheadrightarrow\pi_1(X_{\bar v})$ is not injective, as desired. 
\qed\enddemo

\demo{Proof of Theorem 3.6} (cf. [S], Theorem 4.4.) 
Theorem 3.6 follows easily from Theorem 3.9. Indeed, there exist 
a connected and reduced $k$-scheme $\Cal S'$ of finite type, 
a finite surjective $k$-morphism $\Cal S'\to \Cal S$ and 
a proper and smooth 
%(relative) 
$\Cal S'$-curve $f:\Cal X'\to \Cal S'$ of 
(constant) genus $g$, 
such that the (coarse) classifying morphism $\Cal S'\to\Cal M_{g,k}$ 
for $f$ 
coincides with the composite of $\Cal S'\to \Cal S\subset\Cal M_{g,k}$ 
(cf. [S], Proof of Proposition B.1). Suppose that $\pi_1$ is constant 
on $\Cal S\subset \Cal M_{g,k}$. Then $\pi_1$ is constant on $\Cal S'$ 
for $\Cal X'\to\Cal S'$. Now, by Theorem 3.9, $\Cal X'\to\Cal S'$ is 
isotrivial, which implies that $\Cal S\subset\Cal M_{g,k}$ consists 
of a single point. This is absurd, as $\dim(\Cal S)>0$. 
\qed\enddemo

Next, we work with the assumptions in Theorems 3.12 and 3.13. 
In particular, 
$k$ is a field of characteristic $p>0$,
$S$ is a connected and reduced $k$-scheme of finite type, 
$f:X\to S$ is a non-isotrivial proper and smooth 
%(relative) 
$S$-curve of genus $g$, 
and $D\subset X$ is a relatively \'etale divisor of degree $n$, 
such that $2-2g-n<0$. 

\demo{Proof of Theorem 3.12} (cf. [S], Theorem 4.8.) 
Theorem 3.12 follows easily from Theorem 3.9. Indeed, suppose that 
$\pi_1^t$ is constant on $S$ for $(X,D)/S$. 
By using the tame version 
of Lemma 4.1 (cf. [St], Proposition 2.3), 
we may construct a tame Galois cover $(X',D')/S'$ of $(X,D)/S$ 
which is geometrically of degree prime to $p$ 
and ramified at every point of $D$, such that 
the genus of (the fibres of) $X'$ is $\geq 2$
(cf. [T3], Theorem (8.1)). 
By construction, 
$\pi_1^t$ is constant on $S'$ for $(X',D')/S'$, which implies 
that 
$\pi_1$ is constant on $S'$ for $X'\to S'$. (This can be proved by 
either considering the tame version of the specialization homomorphism 
(cf. Remark 3.11(i)) 
or resorting to [T2], Theorem (5.2).) Now, by Theorem 3.9, 
$X'\to S'$ is isotrivial, which implies that $(X,D)/S$ is isotrivial, 
as in [T3], Proof of Theorem (8.1). This completes the proof. 
\qed\enddemo

\demo{Proof of Theorem 3.13} 
Theorem 3.13 follows from Theorem 3.12, together with [T1], Corollary 1.5. 
\qed\enddemo

$$\text{References.}$$

\noindent
[BLR]  Bosch, S., L\"utkebohmert, W., Raynaud, M., 
N\'eron models, Ergebnisse der Mathematik und ihrer Grenzegebiete (3), 21, Springer-Verlag, 1990.  

\noindent
[EGA III] Grothendieck, A., \'Elements de g\'eom\'etrie alg\'ebrique, 
\'Etude cohomologique des faisceaux coh\'erents, 
Publications math\'ematiques IHES, 11 and 17, 1961 and 1963.

\noindent
[EGA IV] Grothendieck, A., \'Elements de g\'eom\'etrie alg\'ebrique, 
\'Etude locale des sch\'emas et des morphismes de sch\'emas, 
Publications math\'ematiques IHES, 20, 24, 28 and 32, 1964, 1965, 1966 and 1967.

\noindent
[FJ] Fried, M. D., Jarden, M. F., 
Field arithmetic, 
Ergebnisse der Mathematik und ihrer Grenzgebiete (3), 11, 
Springer-Verlag, 
%\publaddr Berlin
1986. 

%\noindent
%[L] Liu, Q., Algebraic geometry and arithmetic curves, Oxford graduate texts in mathematics 6. 
%Oxford University Press, 2002.

\noindent
[R] Raynaud, M., 
Sections des fibr\'es vectoriels sur une courbe, 
Bull. Soc. Math. France 
110 
%no. 1
103--125 
(1982). 

\noindent
[S] Sa\"\i di, M., On complete families of curves with a given fundamental group in positive characteristic, Manuscripta Math. 118, 425--441 (2005).

\noindent
[SGA1] Grothendieck, A., Rev\^etements \'etales et groupe fondamental, Lecture 
Notes in Math. 224, Springer, Heidelberg, 1971.

\noindent
[Sh] Shafarevich, I., On $p$-extensions, A.M.S translation, Serie 2, Vol.4 (1965).

\noindent
[St] Stix, J., A monodromy criterion for extending curves, 
Internat. Math. Res. Notices 2005, No.29, 1787--1802. 

\noindent
[T1] Tamagawa, A., 
On the fundamental groups of curves over algebraically closed fields of characteristic $>0$, 
Internat. Math. Res. Notices 1999, No.16, 853--873. 

\noindent
[T2] Tamagawa, A.,  On the tame fundamental groups of curves over algebraically closed fields of characteristic $>0$, in Galois groups and fundamental groups,
Math. Sci. Res. Inst.
Publ., 41 (L. Schneps, ed.), Cambridge University Press, 2003, pp. 47--105.

\noindent
[T3] Tamagawa, A., Finiteness of isomorphism classes of curves in positive characteristic with prescribed fundamental groups, J. Algebraic Geometry, 13 (2004) 675--724. 

\bigskip
%\newpage
\noindent
Mohamed Sa\"\i di

\noindent
College of Engineering, Mathematics, and Physical Sciences

\noindent
University of Exeter

\noindent
Harrison Building

\noindent
North Park Road

\noindent
EXETER EX4 4QF 

\noindent
United Kingdom

\noindent
M.Saidi\@exeter.ac.uk

\bigskip
\noindent
Akio Tamagawa

\noindent
Research Institute for Mathematical Sciences

\noindent
Kyoto University

\noindent
KYOTO 606-8502

\noindent
Japan

\noindent
tamagawa\@kurims.kyoto-u.ac.jp
\enddocument